\def\versiondate{19 May 2008}
\input math.macros
\input Ref.macros

\checkdefinedreferencetrue
\continuousfigurenumberingtrue
\theoremcountingtrue
\sectionnumberstrue
\forwardreferencetrue
\citationgenerationtrue
\nobracketcittrue
\hyperstrue
\initialeqmacro

\input \jobname.key
\bibsty{myapalike}

\def\subsectionspace{.1 \vsize}
\def\subsectionskip{\smallskip}
\def\bsubsection#1{\vskip0pt plus \subsectionspace \penalty -250
                   \vskip0pt plus -\subsectionspace \bigskip
                   \noindent {\bf #1.} \subsectionskip}

\def \Ghat {{\widehat G}}

\def\E{{\bf E}}
\def\F{{\cal F}}
\def\P{{\bf P}}

\def\vertex{{\ss V}}
\def\V{\vertex}
\def\edge{{\ss E}}
\def\dedge{{\buildrel \rightarrow \over \edge}}
\def\ed#1{[#1]}  
\def\uedg(#1, #2){\ed{#1, #2}}
\def\oedg(#1){\Seq{#1}}  
\def\fsf{{\ss FSF}}         
\def\wsf{{\ss WSF}}         
\def\ust{{\ss UST}}            
\def\fusf{{\ss FUSF}}          
\def\wusf{{\ss WUSF}}          
\def\HH{{\Bbb H}}           
\def\CHI{\raise2pt\hbox{$\chi$}}
\def\smalloplus{\raise1pt\hbox{$\,\scriptstyle \oplus\;$}}
\def\bd{\partial}
\def\ev#1{{\cal #1}}
\def\F{{\cal F}}    
\def\EC(#1,#2){{\cal C}(#1 \leftrightarrow #2)}   
\def\ER(#1,#2){{\cal R}(#1 \leftrightarrow #2)}   
\def\wireEC(#1,#2){{\cal C}^{\rm W}(#1 \leftrightarrow #2)}  
\def\freeEC(#1,#2){{\cal C}^{\rm F}(#1 \leftrightarrow #2)}  
\def\wireER(#1,#2){{\cal R}^{\rm W}(#1 \leftrightarrow #2)}  
\def\freeER(#1,#2){{\cal R}^{\rm F}(#1 \leftrightarrow #2)}  
\def\wire#1{#1^{\rm W}}   
\def\en{{\cal E}}   
\def\etail#1{#1^-}   
\def\ehead#1{#1^+}   

\def\opn#1{\hbox{\rm #1}} 

\def\diam{\opn{diam}} 
\def\gh{G}  
\def\edges{{\ss E}}
\font\frak=eufm10   
\font\scriptfrak=eufm7
\font\scriptscriptfrak=eufm5
\def\mathfrak#1#2{
\def#1{{\mathchoice%
{{\hbox{\frak #2}}}%
{{\hbox{\frak #2}}}%
{{\hbox{\scriptfrak #2}}}%
{{\hbox{\scriptscriptfrak #2}}}}}}
\mathfrak{\fo}{F}  
\def\prfl{\kappa}   
\def\sumC#1{\pi(#1)} 
\def\divergence{\nabla\cdot}
\def\te{\tilde e}
\def\zalpha{\pi}  
\def\bp{o}   
\def\bd{\partial}   
\def\bdv{\bd_\vertex}   
\def\bdvi{\bd_\vertex^{\rm int}}   
\def\bde{\bd_\edge}    
\def\bdei{\bd_\edge^\infty}    

\def\reverse#1{-#1}   
\def\pc{p_{\rm c}}
\def\Aut{{\rm Aut}}

\def\bfz{{\bf 0}}  

\def\BLPSgip{\ref b.BLPS:gip/}

\def\BLPSusf{\ref b.BLPS:usf/, later referred to as BLPS (\refbyear
{BLPS:usf}),%
\def\BLPSusf{BLPS (\refbyear {BLPS:usf})}}


\ifproofmode \relax \else\head{} {Version of \versiondate}\fi
\vglue20pt

\title{Ends in Uniform Spanning Forests}

\author{Russell Lyons, Benjamin J.~Morris, and Oded Schramm}

\abstract{It has hitherto been known that in a transitive unimodular graph,
each tree in the wired spanning forest has only one end a.s.
We dispense with the assumptions of transitivity and unimodularity,
replacing them with a much broader condition on the isoperimetric profile
that 
requires just slightly more than uniform transience.
}

\bottomII{Primary 60B99. 
 Secondary
60D05, 
20F32
.}
{Spanning trees, Cayley graphs.
}
{Research partially supported by 
NSF grants DMS-0406017
and  
DMS-0707144.}

\bsection {Introduction}{s.intro}

The area of uniform spanning forests has proved to be very fertile.
It has important connections to several areas, such as
random walks, sampling algorithms, domino tilings, electrical networks, and potential
theory. It led to the discovery of the SLE processes, which are a major
theme of contemporary research in planar
stochastic processes.
Although much is known about uniform spanning forests, several important
questions remain open.
We answer some of them here.

Given a finite connected graph, $\gh$, let $\ust(\gh)$ denote the uniform
measure on spanning trees of $\gh$.
If an infinite connected graph $\gh$ is exhausted by a sequence of finite
connected subgraphs $\gh_n$, then the weak limit of
$\Seq{\ust(\gh_n)}$ exists.
(This was conjectured by R.~Lyons and proved by
\ref b.Pemantle:ust/.)
However, it may happen that the limit measure is not supported on trees,
but only on forests.
This limit measure is now called the {\bf free (uniform) spanning forest} on
$\gh$, denoted $\fsf$ or $\fusf$.
If $\gh$ is itself a tree, then this measure is
trivial, namely, it is concentrated on $\{\gh\}$. Therefore, \ref
b.Hag:ustmst/ formally
introduced another limit that had been considered on $\Z^d$
more implicitly by \ref b.Pemantle:ust/ and explicitly by
\ref b.Hag:rcust/, namely, the weak limit of the uniform
spanning tree measures on $\gh_n^*$, where $\gh_n^*$ is the graph $\gh_n$ with
its boundary identified (``wired")
to a single vertex. As \ref b.Pemantle:ust/ showed,
this limit also exists
on any graph and is now called the {\bf wired (uniform) spanning forest},
denoted $\wsf$ or $\wusf$.
It is clear that both $\fsf$ and $\wsf$ are concentrated on the set of
spanning forests\Fnote{By a
``spanning forest'', we mean a subgraph without cycles that contains every
vertex.} of $\gh$ all of whose
trees are infinite.
Both $\fsf$ and $\wsf$ are important in their own right; see \ref
b.Lyons:bird/ for a survey and \BLPSusf\ for a comprehensive
treatment.

A very basic global topological invariant of a tree is the number of
ends it has. The ends of a tree can be defined as equivalence 
classes of infinite simple paths in the tree, where two paths are equivalent 
if their symmetric difference is finite.
Trees that have a single
end are infinite trees that are in some sense very close to being finite.
For example, a tree with one end is recurrent for simple random walk
and has critical percolation probability $\pc=1$.
In an attempt to better understand the properties of the \wsf\ and \fsf,
it is therefore a natural problem to determine the number of ends
in their trees.
In fact, this was one of the very first questions asked about infinite
uniform spanning forests.

\ref b.Pemantle:ust/ proved that a.s.,
each tree of the \wsf\ has only one end in 
$\Z^d$ with $d = 3, 4$ and also showed that there are at most two
ends per tree for $d \ge 5$.
\BLPSusf\ completed and extended this to all unimodular transitive
networks, showing that each tree has only one end.
In fact, each tree of the \wsf\ has only one end in every quasi-transitive
transient network, as well as in a host of other natural networks, as we
show in Theorems~\briefref t.profile1end/ and~\briefref t.wsfoneend/ below.
Our proof is simpler even for $\Z^d$, besides having the
advantage of greater generality.
Instead of transitivity, the present proof 
is based on a form of uniform transience arising from an
isoperimetric profile.
The proof relies heavily on electrical networks and not at all on random
walks, in contrast to the previous proofs of results on the number of ends in
the $\wsf$.
We use the fact discovered by \ref b.Morris:wsf/ that with an appropriate
setup, the conductance to infinity is a martingale with respect
to a filtration that examines edges sequentially and discovers if they
are in the tree of the origin.
Our results answer positively Questions 15.3 and 15.5 of \BLPSusf.

The statement that each tree has one end is a qualitative statement.
However, our method of proof can provide quantitative versions of this
result. To illustrate this point, we show that in $\Z^d$, the
Euclidean
diameter of the past of the origin satisfies the tail estimate
$\P[\diam>t] \le C_d\,t^{-\beta_d}$,
where $\beta_d={1\over 2}-{1\over d}$ and the {\bf past} of a vertex $x$
is the union of the finite connected components of $T_x\backslash x$,
where $T_x$ is the $\wsf$ tree containing $x$.
However, this tail estimate is not optimal in $\Z^d$.

For an application of the property of one end to the abelian sandpile model,
see \ref b.JaraiRedig/; more precisely, we prove in \ref l.bi-rooted/
a property they use that is equivalent to one end.
For an illustration of the usefulness of the one-end property
in the context of the euclidean minimal spanning tree, see \ref
b.Krikun:alloc/.

\bsection{Background, Notations and Terminology}{s.back}

We shall now introduce some notations and give a brief background on
electrical networks and uniform spanning forests.
For a more comprehensive account of the background, please see \BLPSusf. 

\bsubsection {Graphs and networks}

A {\bf network} is a pair $(G,c)$, where $G=(\vertex,\edge)$ is a graph and
$c : \edge \to (0, \infty)$ is a positive function, which is often called the weight function or the
edge conductance.
If no conductance is specified for a graph, then we take $c \equiv 1$ as
the default conductance; thus, any graph is also a network.
Let $G=(\vertex,\edge,c)$ be a network,
and let $\dedge=\bigl\{\oedg(x, y) \st \uedg(x, y)\in\edge\bigr\}$
denote the set of oriented edges.
For $e\in\dedge$, let $\reverse e$ denote its reversal,
let $\etail e$ denote its tail, and let $\ehead e$ denote its head.

For a set of vertices $K\subset\vertex$, let $\bde K$ be its {\bf edge boundary} that consists
of edges exactly one of whose endpoints is in $K$.
Sometimes $K$ is a subset of two graphs; when we need to indicate in which
graph $G$ we take the edge boundary, we write $\bd_{\edge(G)}(K)$.
Let $G\backslash K$ denote the graph $G$ with
the vertices $K$ and all the edges incident with them removed,
and let $G/K$ denote the graph $G$ with the vertices in $K$ identified (wired) to a
single point and any resulting loops (edges with $\ehead e=\etail e$) 
dropped.
(Note that it may happen that $G/K$ is a multigraph even if $G$ is a simple graph;
that is, $G/K$ may contain multiple edges joining a pair of vertices even
if $G$ does not.)
If $H$ is a subgraph of $G$, then $G\backslash H$ and $G/H$ mean
the same as $G\backslash \vertex(H)$ and $G/\vertex(H)$,
respectively.  However, when $e$ is an edge (or a set of edges),
$G\backslash e$ means $G$ with $e$ removed but no vertices
removed.

Write $|F|_c := \sum_{e\in F} c(e)$ for any set of edges $F$.
If $K$ is a set of vertices, we also write
$$
\zalpha(K):=\sum\Bigl\{c(e)\st e\in\dedge,\, \etail e\in K\Bigr\}.
$$
We shall generally consider only networks where $\zalpha(x)<\infty$ for every
$x\in\vertex$.

On occasion, we shall need to prove statements about infinite networks from
corresponding statements about finite networks. For this purpose, the
concept of an exhaustion is useful. An {\bf exhaustion} of $G$
is an increasing sequence of finite connected subnetworks
$H_0\subset H_1\subset H_2\subset\cdots$
such that $\bigcup_{n\ge 1} H_n=G$.

Given two graphs $G = (\vertex, \edge)$ and $G' = (\vertex', \edge')$, call
a function $\phi : \vertex \to \vertex'$ a {\bf rough isometry} if there are
positive constants $\alpha$ and $\beta$ such that for all $x, y\in \vertex$,
$$
\alpha^{-1} d(x, y) - \beta \le d'\big(\phi(x), \phi(y)\big)
\le \alpha d(x, y) + \beta
\label e.roughisom
$$
and such that every vertex in $G'$ is within distance $\beta$ of the image of
$\vertex$. Here, $d$ and $d'$ denote the usual graph distances on $G$ and
$G'$. The same definition applies to metric spaces, with ``vertex"
replaced by ``point".

\bsubsection {Effective resistance and effective conductance}

Let $(G, c)$ be a connected network.
The {\bf resistance} of an edge $e$ is $1/c(e)$ and will be denoted by
$r(e)$.
A function $\theta:\dedge\to\R$ is called {\bf antisymmetric}
if $\theta(e)=-\theta(\reverse e)$ for all $e\in\dedge$.
For antisymmetric functions $\theta,\theta':\dedge\to\R$, set
$$
(\theta,\theta')_r := {1\over 2} \sum_{e\in \dedge}
r(e)\,\theta (e)\,\theta'(e) 
\,.
$$
The {\bf energy} of $\theta$ is given by $\en(\theta) :=(\theta,\theta)_r$.
The {\bf divergence} of $\theta$ is
the function $\divergence\theta:\vertex\to\R$ defined by
$$
\divergence \theta(x):=
\sum\Bigl\{ \theta(e) \st e\in\dedge,\,\etail e = x\Bigr\}.
$$
(In some papers the definition of divergence differs by a factor of $\sumC x ^{-1}$
from our current definition.) 

If $f:\vertex\to\R$ is any function, its {\bf gradient}
$\nabla f$ is the antisymmetric function on $\dedge$
defined by $\nabla f(e):=c(e)\,\bigl(f(\ehead e)-f(\etail e)\bigr)$.
The {\bf Dirichlet energy} of $f$ is defined as
$D(f):=(\nabla f,\nabla f)_r$, and its {\bf laplacian} is 
$\Delta f:=\divergence\nabla f$. The function $f$ is {\bf harmonic}
on a set $A\subset\vertex$ if $\Delta f=0$ on $A$.

Consider now the case where $G$ is finite.
Let $f:\vertex\to\R$ and let $\theta:\dedge\to\R$ be antisymmetric.
The useful identity
$$
(\nabla f,\theta)_r = -\sum_{v\in\vertex} f(v)\,\divergence \theta(v)
\label e.ibp
$$
follows by gathering together the terms involving $f(v)$ 
in the definition of $(\nabla f,\theta)_r$.
Let $A$ and $B$ be nonempty disjoint subsets of $\vertex=\vertex(G)$.
A unit {\bf flow} from $A$ to $B$ is an antisymmetric
$\theta:\dedge\to\R$ satisfying $\divergence \theta(x)=0$
when $x\notin A\cup B$, and $\sum_{x\in A}\divergence\theta(x)=1$.
Since, clearly, $\sum_{x\in\vertex} \divergence\theta(x)=0$,
it follows that $\sum_{x\in B}\divergence\theta(x)=-1$.
The {\bf effective resistance} between $A$ and $B$ in $(G,c)$
is the infimum of $(\theta,\theta)_r$ over all unit flows $\theta$ from
$A$ to $B$, and will be denoted by $\ER(A,B)$.
This minimum is achieved by the {\bf unit current flow}.
The {\bf effective conductance} between $A$ and $B$,
denoted $\EC(A,B)$, is the infimum of $D(f)$ over all functions $f:\vertex\to\R$
satisfying $f=0$ on $A$ and $f=1$ on $B$.
It is well known (and follows from~\ref e.ibp/) that $\EC(A, B) = \ER(A, B)^{-1}$.
Furthermore, $D(f)$ is minimized for a function $f$ that is harmonic except
on $A \cup B$ and whose gradient is proportional to the unit current flow.
\comment{ 
For the sake of readers
who are unfamiliar with this, we give the proof.
Suppose that $\theta$ is a unit flow from $A$ to $B$
and $f:\vertex\to\R$ satisfies $f=0$ on $A$ and $f=1$ on $B$.
Then by~\ref e.ibp/ we have 
$$
(\theta,\nabla f)_r=1\,,
\label e.ip
$$
since $\divergence\theta=0$ except on $B\cup A$,
and $\sum_{x\in B}\divergence \theta(x)=-1$.
Hence, Cauchy-Schwarz gives $\en(\theta)\,D(f)\ge 1$.
Let $h$ be the function that is $0$ on $A$, $1$ on $B$,
and harmonic in $\vertex\setminus(A\cup B)$.
(The function $h$ has a probabilistic description: $h(x)$ is
the probability to hit $B$ before $A$ for the random walk starting at $x$
with transition probabilities at $y$ proportional to $c(e)$ for $e$
incident to $y$.)
Then $\Delta h=\divergence\nabla h$ is $0$ on
$\vertex\setminus (A\cup B)$. If $a:=\sum_{x\in A}\Delta h$,
then $a^{-1}\,\nabla h$ is a unit flow from $A$ to $B$.
Hence,~\ref e.ip/ gives $(a^{-1}\,\nabla h,\nabla h)_r=1$,
which implies $D(h)=a$ and $\en(a^{-1}\,\nabla h)=a^{-1}$.
Since $\en(\theta)\,D(f)\ge 1$ holds for every $\theta$ and $f$
as above, we conclude that $h$ minimizes $D(f)$ among such $f$, that
$a^{-1}\,\nabla h$ minimizes $\en(\theta)$ among such $\theta$
and that $\EC(A,B)=\ER(A,B)^{-1}$.
}

Now suppose that the network $(G,c)$ is infinite.
The effective conductance from a finite set $A\subset\vertex$ to $\infty$
is defined as the infimum of $D(f)$ over all $f:\vertex\to \R$
such that $f=0$ on $A$ and $f=1$ except on finitely many vertices.
This will be denoted by $\EC(A,\infty)$, naturally.
The effective resistance to $\infty$ can be defined
as $\ER(A,\infty):=\EC(A,\infty)^{-1}$, or, equivalently,
as the infimum energy of any unit flow from $A$ to $\infty$,
where a unit flow from $A$ to $\infty$ is an antisymmetric 
$\theta:\dedge\to\R$ that satisfies $\divergence \theta=0$ outside of
$A$ and $\sum_{x\in A}\divergence\theta(x)=1$.
When $B\subset\vertex$, we define
$\ER(A,B\cup\{\infty\})$ as the infimum energy of any antisymmetric
$\theta:\dedge\to\R$ whose divergence vanishes in $\vertex\setminus(A\cup B)$
and which satisfies $\sum_{x\in A}\divergence\theta(x)=1$,
and define $\EC(A,B\cup\{\infty\})$ as
$\ER(A,B\cup\{\infty\})^{-1}$, or, equivalently,
as the infimum of $D(f)$, where $f$ ranges over all
functions that are $0$ on $A$, $1$ on $B$, and different from $1$
on a finite set of vertices.

\comment{
We remark that there are two natural definitions for
the effective conductance between $A$ and $B$ in an
infinite network $G$: the free
effective conductance is the limit of $\EC(A,B;H_n)$,
where $H_n$ is an exhaustion of $G$,
and the wired effective conductance is  the limit
of $\EC(A,B;G/(G\setminus H_n))$. It
is not hard to see, using monotonicity, that the limits
exist and do not depend on the choice of the exhaustion;
however, in general the limits may be different.
Similarly, there is a wired and a free definition
of the effective resistance in an infinite network,
which may or may not agree.
}

When $A$ or $B$ belong to two different networks
under consideration, we use the
notations $\ER(A,B;G)$ and $\EC(A,B;G)$ in order to specify 
that the the effective resistance or conductance is with respect
to the network $G$.

\bsubsection {Uniform spanning trees and forests}

If $(G,c)$ is finite,
the corresponding {\bf uniform spanning tree} is the measure on spanning
trees of $G$ such that the probability of a spanning tree $T$ is
proportional to $\prod_{e \in T} c(e)$.

The following relation between spanning trees and
electrical networks is due to \ref b.Kirchhoff/.

\procl p.einT  Let $T$ be a uniform spanning
tree of a finite network $G$ and $e$ an edge of $G$.  Then
$$
\P[e\in T] 
= i(e) = c(e)\,\ER(\etail{e},\ehead{e})
\,,
$$
where $i$ is the unit current from $\etail e$ to $\ehead e$.
\endprocl

If $e$ is an edge in $G$, it is easy to see that
the conditional law of the uniform spanning tree $T$ given $e\in T$ is the same
(considered as a set of edges) as the law of the uniform spanning
tree of $G/e$ union with $\{e\}$.
Also, the conditional law of $T$ given $e\notin T$ is
the same as the law of the uniform spanning tree
of $G\backslash e$.
These facts will be very useful for us.

Now assume that $(G,c)$ is infinite,
and let $H_n$ be an exhaustion of $G$.
Let $H_n^*$ denote the graph $G/(G\backslash H_n)$,
namely, $G$ with the complement of $H_n$ identified
to a single vertex (and loops dropped).
Let $T_n$ denote the uniform spanning tree on
the network $(H_n,c)$, and let $T_n^*$ denote
the uniform spanning tree on $(H_n^*,c)$.
Monotonicity of effective resistance and~\ref p.einT/ imply that for
every $e\in\edge$ the limit of $\P[e\in T_n]$
exists as $n\to\infty$.  By the previous paragraph,
it is easy to conclude that the weak limit of the
law of $T_n$ exists (as a measure on the Borel
subsets of $2^\edge$). This measure is
the {\bf free uniform spanning forest} on $G$,
and is denoted by $\fsf$.
Likewise, the weak limit of the law of $T_n^*$
exists (this time the monotonicity goes in the
opposite direction, though), and is
called the {\bf wired uniform spanning forest}
on $G$, which is denote by $\wsf$.
In $\Z^d$ and many other networks, we have
$\wsf=\fsf$. However, there are some interesting
cases where $\wsf\ne\fsf$
(an example is when $G$ is a transient tree
and $c\equiv 1$).

In the sequel, stochastic domination of probability measures on spanning
forests will refer to the partial order induced by inclusion when forests
are regarded as sets of edges.

\comment{ 
Samples from the $\wsf$, the $\fsf$ and from the
uniform spanning tree in a finite network will
generally be denoted by $\fo$ below.
}

\bsection{$\wsf_o$}{s.rootwired}

In this section, $G$ is an arbitrary connected network and $o$
is some vertex in $G$.
We now define a probability measure on spanning forests
of $G$, which is the wired spanning forest with $o$ wired to $\infty$.
Suppose that $G$ is exhausted by finite
subgraphs $\Seq{G_n}$. Let $\Ghat_n$ be the graph obtained from 
$G$ by identifying $o$ and the exterior of $G_n$ to 
a single point.
Then the wired spanning forest on $G$ with $o$ wired to
$\infty$ is defined as the weak limit (as a set of edges) of the 
uniform spanning tree on $\Ghat_n$
as $n\to\infty$,
and will be denoted by $\wsf_o$.
The existence of the limit follows from monotonicity
by the same argument that gives the existence of the $\wsf$.

\procl p.vertex-wired
Let $\gh$ be a transient network and $\bp \in \gh$.
Given a forest $\fo$ in $\gh$, let $\fo(\bp)$ be its component that contains
$\bp$.
Then $|\fo(\bp)| < \infty$ $\wsf_\bp$-a.s.\ iff $\wsf$-a.s., there do not exist two
edge-disjoint infinite paths in $\fo(\bp)$ starting at $\bp$.
\endprocl

The ``if" direction of this proposition was proved in \ref b.JaraiRedig/.
The ``only if" direction
(which is the direction we shall use)
is an immediate consequence of the following
lemma in which we take $x$ to be $\bp$ and $y$ to be the 
 wired vertex in
$G/\bigl(\vertex(G\backslash G_n)\backslash\{\bp\}\bigr)$,
and then take the weak limits:

\procl l.bi-rooted
Let $\gh$ be a finite connected network and $x, y \in \vertex$ be distinct
vertices.
Given a spanning tree $T$ of $\gh$, let $L(T)$ be the path in $T$ that
connects $x$ to $y$.
The uniform spanning tree in $\gh/\{x, y\}$ stochastically dominates $T
\backslash L(T)$ when $T$ is a uniform spanning tree of $\gh$.
\endprocl

Before presenting the proof, we recall that for any tree $T_0\subset \gh$,
the set of edges
in the uniform spanning tree $T$ conditioned on $T_0\subset T$ has the same distribution as the union
of $\edges(T_0)$ with the set of edges of the uniform spanning tree of 
$G/\vertex(T_0)$.
Here, $\vertex(T_0)$ and $\edges(T_0)$ denote the vertices 
and edges of $T_0$, respectively.

\proof
Indeed, with the proper identification of edges,
conditioned on $L(T)$,
$T \backslash L(T)$ is the uniform spanning tree of
$\gh/\vertex\big(L(T)\big)$.
Since $\vertex\big(L(T)\big)$ includes $x$ and $y$, it follows from
a repeated application of the well-known 
negative association theorem of \ref b.FedMih/ (see, e.g., Theorem 4.4
of \BLPSusf) that the stochastic domination holds when we condition on $L(T)$.
By averaging we conclude that it also holds unconditionally.
\Qed

\comment{
\procl l.dominate
If $\wsf_\bp$-a.s.\ the component of $\bp$ in $\fo$ is finite,
then $\wsf$-a.s., there do not exist two
edge-disjoint infinite paths in $\fo$ starting at $\bp$.
\endprocl

\ref b.JaraiRedig/ also considered $\wsf_o$, and proved the converse of
\ref l.bi-rooted/.

In the proof of the lemma, we shall use the following monotonicity result
of \ref b.FedMih/ (see also Theorem 4.4 of \BLPSusf).
Let $\gh$ be a finite connected graph and let ${\cal F}$ be some fixed
collection of
subsets of $\edges(\gh)$. Suppose
${\cal F}$ is monotone increasing in the sense that
whenever $F_1\subset F_2\subset\edges(\gh)$
and $F_1\in {\cal F}$, we also have $F_2\in{\cal F}$.
Let $T$ be the set of edges in the uniform spanning tree of $\gh$, and let
$S\subset\edges(\gh)$ be a set of edges such that $\P[S\subset T]>0$.
Then $\P[T\backslash S\in{\cal F}\mid S\subset T]\le \P[T\backslash S\in{\cal F}]
\le \P[T\in{\cal F}]$.

\proof
Let $H$ be a connected finite subgraph of $G$ containing $o$,
and let $H_o:=H\backslash \{o\}$.
Let $H^*$ and $H_o^*$ denote the wired graphs $G/(G\backslash H)$
and $G/(G\backslash H_o)$, respectively. In each case, we denote
the wired vertex by $\infty$.
Let $\fo$ denote the uniform spanning tree on $H^*$,
let $\fo_o$ denote the uniform spanning tree on $H^*_o$,
and let $\gamma$ denote the simple path in $\fo$
joining $o$ and $\infty$. 
Conditional on $\gamma$, the edges of the uniform spanning tree of
$H^*$ that are not in $\gamma$ have the same conditional distribution
as the edges of the uniform spanning tree of $H^*/\gamma$
and hence of $H^*_o/\gamma$.
By monotonicity, for each fixed $k$, the conditional probability given $\gamma$
that there is a path of length $k$ in $H$ which has $o$
as an endpoint and whose edges are in $\fo\backslash\gamma$
is no larger than the probability that there is a path of length
$k$ in $H$ which has $o$ as an endpoint and whose edges
are in $\fo_o$.
Therefore, the latter also bounds the probability that there are two edge-disjoint
paths of length $k$ in $\fo$ starting at $o$.
The result follows by first letting $H$ exhaust $\gh$, and then letting
$k\to\infty$.
\Qed
}

Informally, the following lemma says that conductance to infinity is
a martingale.
When $\fo$ is a spanning forest of a graph $G$ and $E$ is a set of edges of
$G$, we denote the graph $\big(\vertex(G), \edges(\fo) \cap E\big)$ by
$\fo \cap E$.

\procl l.martingale
Let $\fo$ be a sample from $\wsf_o$ on $(G, c)$.
Let $E_0\subset E_1$ be finite sets of edges in $G$, and
for $j=0,1$, let $S_j$ be the set of vertices of
the connected component of $o$ in $\fo\cap E_j$,
and let $M^j$ be the effective conductance from $S_j$ to $\infty$
in $G\backslash E_j$.
On the event that every edge in $E_1$ has at least one endpoint in
$S_0$, we have
$$
\E\bigl[M^1\bigm| \fo\cap E_0\bigr]=M^0\,.
$$
\endprocl

\proof
Suppose that $G$ is exhausted by finite subgraphs $\Seq{G_n}$,
with 
$E_1 \subset \edges(G_n)$ for all $n$. Let $\Ghat_n$ be defined as in the 
definition of $\wsf_o$ and 
let $T_n$ be a random spanning tree of 
$\Ghat_n$. 
For $j=0,1$ and $n \geq 1$, 
let $\Ghat^j_n$ be the graph obtained from $\Ghat_n$ by contracting or 
deleting the edges in $E_j$ according to whether they 
are in $T_n$ or not and 
let $G^j$ be the graph obtained from $G$ by contracting or 
deleting the edges in $E_j$ according to whether they 
are in $\fo$ or not;
let $M_n^j$ be the effective conductance from $o$ to the exterior 
of $G_n$ in 
$G^j$, and  
let $\bar M^j$ be the effective conductance
between 
$\bp$ and $\infty$ in $G^j$.  
Theorem 7 of 
\ref b.Morris:wsf/ 
implies that 
$$
\E\bigl[M_n^1 \bigm| \Ghat_n^0\bigr] = M_n^0\,.
\label e.condmart
$$
Since the weak 
limit of the $T_n$ 
is $\wsf_o$, 
taking the limit
of both sides of equation \ref e.condmart/ as $n \to \infty$ 
 gives $\E\bigl[\bar M^1 \bigm| G^0\bigr] = \bar M^0$. 
Note that $\sigma(G^0) = \sigma(\fo \cap E_0)$.

Since on the event that
every edge in $E_1$ has at least one endpoint in $S_0$, 
the graph obtained from $G$ by identifying 
the vertices in $S_j$ and deleting the edges in $E_j$ is 
$G^j$ for $j = 0, 1$, we also have that $\bar M^j = M^j$ on that event,
whence the lemma follows. 
\Qed 
\bsection{The case of $\Z^d$}{s.Zdcase}

This section is devoted to proving the following theorem.

\procl t.tail
Let $d>2$, $d\in\N$, and let $\fo$ denote a sample from the $\wsf$
on $\Z^d$.  Let $\fo(\bfz)$ denote the connected component of $\bfz$ in $\fo$.
Then $\fo(\bfz)$ a.s.\ has one end, and therefore $\fo(\bfz)\backslash\{\bfz\}$ has
just one infinite connected component a.s.
Moreover, if $Q$ denotes the union of the finite connected components
of $\fo(\bfz)\backslash\{\bfz\}$,
then for all $t>0$,
$$
\P\bigl[\diam(Q)>t\bigr] 
\le C_d\,t^{-\beta_d},
$$
where $\beta_d:={1\over 2}-{1\over d}$ and $C_d$ is some constant depending only on $d$.
Here, $\diam$ means the Euclidean diameter.
\endprocl

Of course, the first statement, namely that $\fo(\bfz)$ has one end,
is not new, but the proof we give is rather different from the proof
in \BLPSusf.
Our proof, without the quantitative estimate, is significantly shorter than
that in \BLPSusf.

The tail estimate on $\diam(Q)$ that the theorem provides
is not optimal.
It is possible to show that $\P[\diam(Q)>t]$ behaves like $t^{-2}$ 
in $\Z^d$, $d>4$, with
possible polylogarithmic corrections, but we shall not prove this in the present paper.
However, the proof we have in mind for the stronger estimates relies heavily on the fact
that $\Z^d$ is transitive, and would not work for non-transitive networks
that are similar to $\Z^d$,
such as $\N\times \Z^{d-1}$ or $\Z^d$ with variable edge conductances bounded between
two positive constants. In contrast, it is not too hard to see
that the proof of \ref t.tail/ given below easily extends to these settings.

\medskip

Throughout this section, $O_d(1)$ will stand for an unspecified positive
constant depending only on the dimension $d$.

\procl l.quantitativeS
Let $d\ge 3$ and let $G$ be the graph $\N\times\Z^{d-1}\subset\Z^d$.
If $S\subset\vertex(G)$ is finite
and nonempty and contained in $\{0\}\times \Z^{d-1}$,
then the effective conductance from
$S$ to $\infty$ in $G$ is at least 
$|S|^{{d-2\over d-1}}/O_d(1)$.
\endprocl

\proof
We consider first $\EC(S,\infty;\Z^d)$.
For $v\in\Z^d$, let $g_v$ be the Green function in $\Z^d$, that is,
$g_v(x)$ is the expected number of visits to $x$ for simple random walk
started at $v$.
Set $\theta_v:=-\nabla g_v/(2 d)$, that is, $\theta_v(e)=
\big(g_v(\etail e)-g_v(\ehead e)\big)/(2 d)$
for oriented edges $e$.
It is easy to verify that the divergence of $\theta_v$ is zero
at every $u\in \Z^d\setminus\{v\}$ and 
that $g_v(v)=1+g_u(v)=1+g_v(u)$ holds whenever $u$ neighbors $v$.
Therefore,
$\theta_v$ is a unit flow from $v$ to $\infty$ in $\Z^d$.
If we formally apply~\ref e.ibp/ with
$f:=g_u$ and $\theta:=\theta_v$, we get $(\theta_v,\theta_u)_1=g_u(v)/(2 d)$.
However, since in this case the summation corresponding to~\ref e.ibp/ is
infinite, we need to be more careful.
Let $G$ be exhausted by $\Seq{G_n}$ and
$\wire{G_n} := G/(G\backslash G_n)$.
We may assume that $v \in \vertex(G_n)$ for all $n$.
Let $g_{v, n}$ be the Green function in $\wire{G_n}$
for the walk on $\wire{G_n}$ killed when it first leaves $G_n$, so that it
takes the value 0 off of $G_n$.
Clearly $\lim_{n \to\infty} g_{v, n} = g_v$ pointwise.
Also, $g_{v, n}$ is harmonic on $\vertex(G_n) \setminus \{v\}$.
Thus, $\theta_{v, n} = - \nabla g_{v, n}/(2 d)$ is the unit current
flow on $\wire{G_n}$ from $v$ to the complement of $G_n$.
Because $g_{v, n} \to g_v$, we have $\theta_{v, n} \to \theta_v$ on each
edge.
Fatou's lemma implies that $\en(\theta_v) \le \liminf_{n \to\infty}
\en\big(\theta_{v, n}\big)$. 
If $\en(\theta_v) < \limsup_{n \to\infty}
\en\big(\theta_{v, n}\big)$, 
then for some $n$, the restriction of $\theta_v$
to $\edges(\wire{G_n})$ would give a unit flow from $v$ to the complement
of $G_n$ with smaller energy than $\theta_{v, n}$, a contradiction, whence
$\en(\theta_v) = \lim_{n \to\infty} \en\big(\theta_{v, n}\big)$.
It follows that $\en(\theta_{v, n} - \theta_v) \to 0$ as $n \to\infty$.
Since $(\theta_{v, n},\theta_{u, n})_1=g_{u, n}(v)/(2 d)$,
the result follows by taking limits.
\comment{
This is not too difficult. We can pick an $\epsilon>0$ very small,
and let $H_\epsilon$ be the subgraph of $\Z^d$ spanned by the
vertices $w$ satisfying $g_v(w)\ge\epsilon$. The graph $H_\epsilon$ is finite,
and we may therefore apply~\ref e.ibp/ to the restrictions
of $\nabla g_v$ and $g_u$ to $H_\epsilon$. Let $U_\epsilon$ denote the vertices
of $H_\epsilon$ that in $\Z^d$ have a neighbor outside of $H_\epsilon$.
Then the restriction $\tilde\theta$ of $\nabla g_v$ to $H_\epsilon$
satisfies $\divergence\tilde\theta\ge 0$ on $U_\epsilon$
and $\sum_{x\in U_\epsilon}\divergence\tilde\theta(x)=1$.
Since $g_u(x)\to 0$ as $|x|\to\infty$, it follows that
$$
\lim_{\epsilon\to0}\sum_{x\in U_\epsilon} \divergence\tilde\theta(x)\,g_u(x)=0\,.
$$
Therefore, when applying~\ref e.ibp/ on $H_\epsilon$ to $\tilde\theta$ and $g_u$,
we get in the limit as $\epsilon\to 0$
$$
(\theta_v,\theta_u)=g_u(v)\,.
$$
}

Define $\theta:=|S|^{-1}\sum_{v\in S}\theta_v$.
Then $\theta$ is a unit flow in $\Z^d$ from $S$ to $\infty$.
Consequently,
$$
\EC(S,\infty; \Z^d)\ge (\theta,\theta)_1^{-1}=
\Bigl( |S|^{-2}\,\sum_{v,u\in S} (\theta_v,\theta_u)_1\Bigr)^{-1}
=
2 d |S|^2\,\Bigl(\sum_{v,u\in S}g_v(u)\Bigr)^{-1}
.
$$
Recall that $g_v(x)\le O_d(1) \, |v-x|^{2-d}$.
(See, e.g., Theorem 1.5.4 of \ref b.Lawler:book-inter/.)
Since $S\subset \{0\}\times\Z^{d-1}$,
it follows that for every $v\in\Z^d$ we
have
$$
\sum_{u\in S} g_v(u)\le O_d(1)\,|S|^{1/(d-1)}
$$
since our upper
bound on $g_v$ is monotone decreasing in distance from $v$, whence the sum
of the bounds is maximized by the set $S$ closest to $v$.
This gives $\EC(S,\infty;\Z^d)\ge |S|^{(d-2)/(d-1)}/O_d(1)$.

The corresponding result now follows for $G$, since we may restrict
$\theta$ to $\edge(G)$ and double it on edges that are not
contained in the hyperplane $\{0\}\times\Z^{d-1}$. The result
is a unit flow from $S$ to $\infty$ in $G$, and its energy is at most
$4$ times the energy of $\theta$ in $\Z^d$. The lemma follows.
\Qed

As a warm up, we start with a non-quantitative version of \ref t.tail/,
proving that a.s.\ the connected components of the $\wsf$ in $\Z^d$
($d\ge 3$) all have one end. 
Let $\fo$ be a sample from $\wsf_\bfz$.
By \ref p.vertex-wired/,
it suffices to show that $\wsf_\bfz$-a.s.\ the connected component of $\bfz$ in
$\fo$ is finite.

Set $B_r:=\{z\in\R^d\st \|z\|_\infty\le r\}$.
We inductively construct an increasing sequence 
$E_0\subset E_1\subset E_2\subset \cdots$ of sets of edges.
Put $E_0:=\emptyset$. Assuming that $E_n$ has been defined, we
let $S_n$ be the vertices of the connected component of $\bfz$ in $\fo\cap E_n$.
If all the edges of $\Z^d$ incident with $S_n$ are in $E_n$,
then set $E_{n+1}:=E_n$.
Otherwise, we choose some edge $e\notin E_n$ incident with $S_n$
and set $E_{n+1}:=E_n\cup\{e\}$.
Among different possible choices for $e$, we take one that minimizes
$\min\{r\st e\subset B_r\}$,
with ties broken in some arbitrary but fixed manner.
Let $M_n:=\EC(S_n,\infty;\Z^d\setminus E_n)$ denote the effective conductance from $S_n$ to
$\infty$ in $\Z^d \backslash E_n$.
Let $\F_n$ denote the $\sigma$-field generated by
$E_n$ and $E_n\cap\fo$.
By \ref l.martingale/, $M_n$ is a martingale
with respect to the filtration $\F_n$,
that is, $\E\bigl[M_{n+1}\bigm| \F_n\bigr]=M_n$.
Since $M_n\ge0$, it follows that $\sup_n M_n<\infty$ a.s.

Let $\ev A_n$ be the event that $E_{n+1}=E_n$, and
set $\ev A:= \bigcup_{n\in\N} \ev A_n$, which is the
event that the component of $\bfz$ in $\fo $ is finite.
For $r \ge 1$, let $n_r$ denote the first $n$ such that $E_n$ contains all the
edges inside $B_r$ that are incident with $S_n$.
(Thus we have $E_{n_r}\subset B_r$ and
either $E_{n_r+1}=E_{n_r}$ or $E_{n_r+1}\not\subset B_r$.)
We claim that for every $c>0$ there is a $\delta_c>0$ such that
for every $r\in\N$ 
$$
\P\bigl[\ev A_{n_{r+1}}\bigm| \F_{n_r} \bigr]\ge \delta_c\,\I{\{M_{n_r}\le c\}}.
\label e.maydie
$$
To prove this, fix some $c>0$ and $r \in \N$,
set $m:=n_r$,  and suppose that $M_m\le c$.
Let $\ell$ be the number of edges 
that connect $S_m$ to a vertex in $[r+1,\infty)\times \Z^{d-1}$,
and let $k_r$ be the number of edges that connect $S_m$ to a vertex
outside $B_r$.
Because of the assumption $M_m \le c$,
\ref l.quantitativeS/ gives a finite upper bound
$K_{c}$ for $\ell$, which depends only on $c$ and $d$.
By symmetry, it follows that 
$k_r\le 2\,d\, K_c$.
If $k_r = 0$, then $\ev A_{n_{r+1}}$ occurs, so \ref e.maydie/ certainly
holds for every $\delta_c \le 1$. Otherwise, suppose that $k_r \ge 1$.
Fix some $j\in\N\cap[m+1,m+k_r]$.
If it so happens that $E_{j-1}\cap\fo\subset E_m$,
then the edge $e_j$ in $E_{j}\setminus E_{j-1}$
is one of those $k_r$ edges that connect $S_m$ to
the complement of $B_r$. Suppose that this
is the case, and let $v$ be the endpoint of $e_j$
that is not in $S_m$.
There is a universal lower bound $\alpha>0$ 
on the conductance from $v$ to $\infty$ in $\Z^d\backslash S_m$.
Thus, the effective conductance between $v$ and
$S_m\cup\{\infty\}$ in the complement of $E_{j-1}$ is at least $1+\alpha$
(seen, e.g., by minimizing Dirichlet energy).
This gives
$$
\wsf_\bfz\bigl[ e_j\in \fo \bigm| \F_{j-1}  \bigr]\,\I{E_{j-1}\cap \fo\subset E_m}
\le (1+\alpha)^{-1}
$$
by \ref p.einT/.
Induction on $j$ yields
$$
 \wsf_\bfz\bigl[ E_j\cap \fo\subset E_m\bigm| \F_m \bigr] \ge \bigl(1- (1+\alpha)^{-1}\bigr)^j
$$
for $j\in\N\cap[m+1,m+k_r]$.
Since $\ev A_{n_{r+1}}$ is the event $E_{m+k_r}\cap\fo\subset E_m$
and since $k_r\le 2\,d\,K_c$, we therefore get \ref e.maydie/ with 
$$
\delta_c = \bigl(1-(1+\alpha)^{-1}\bigr)^{ 2\,d\,K_c}.
$$

Induction on $r$ and \ref e.maydie/ give
$$
\P\Bigl[\sup_{n\le n_r} M_n\le c,\,\neg \ev A_{n_r}\Bigr] \le
(1-\delta_c)^r.
$$
Hence, $\P\bigl[\sup_n M_n\le c,\,\neg\ev A\bigr]=0$.
Because $\sup_nM_n<\infty$ a.s., this clearly implies
that $\wsf_\bfz[\ev A]=1$, which proves that all
the connected components of the $\wsf$ in $\Z^d$ have one end a.s.

\medskip

The above proof can easily be made quantitative, but the bound it provides is
rather weak.  We now proceed to establish a more reasonable bound.

\proofof t.tail
Here, we shall again use the sequence $E_n$ constructed above, as well as the
notations $S_n, M_n, B_r$, etc.
However, for the following argument to work, we need to be more specific
about the way in which $E_{n+1}$ is chosen given $E_n$ and $\fo\cap E_n$.
When $n=n_r$ for some $r$ and the set of edges adjacent to $S_n$ outside of
$B_r$ is nonempty, we now require that $E_{n+1}=E_n\cup\{\te_r\}$,
where $\te_r$ is an edge along which the unit current flow
from $S_n$ to $\infty$ in the complement of $E_n$ is maximal.

Fix some $r\in\N$ and let $m:=n_r$.
Suppose that $E_{m+1}=E_m\cup\{\te_r\}$,
that is, $\ev A_{n_r}$ does not hold.
Let $\theta$ denote the unit current flow from $S_m$ to $\infty$ in the
complement of $E_m$.
Then $M_m=\en(\theta)^{-1}$.
Clearly, $|\theta(\te_r)|\ge 1/k_r$
(where $k_r$ is the number of edges not in $E_m$ that
are incident to $S_m$).
If we have $\te_r\in\fo$, then the restriction of $\theta$ to the
complement of $\te_r$ is a unit flow from $S_{m+1}$ to $\infty$
in the complement of $E_{m+1}$. Therefore, in this case,
$$
M_{m+1}\ge \bigl(\en(\theta)-\theta(\te_r)^2\bigr)^{-1}
\ge \bigl(M_m^{-1}- k_r^{-2}\bigr)^{-1}
\ge M_m\,\Bigl(1+{M_m\over k_r^2}\Bigr).
$$
We also know that $k_r\le O_d(1)\,M_m^{{d-1\over d-2}}$
from \ref l.quantitativeS/.
Consequently, on the event $\te_r\in \fo$, we have
$$
M_{m+1}\ge M_m+M_m^{-2/(d-2)}/O_d(1)\,.
\label e.whener
$$
By \ref p.einT/,
$$
\P\bigl[ \te_r\in\fo \bigm| \F_{n_r}\bigr] \ge (2\,d)^{-1}
\qquad \hbox{ on the event } \neg\ev A_{n_r}\,,
\label e.per
$$
since the conductance of $\te_r$ is $1$ and the effective conductance
in $\Z^d$ 
between the endpoint $x$ of $\te_r$ outside of $S_m$ and $S_m\cup\{\infty\}$
is at most $(2\,d)^{-1}$, as there are $2\,d$ edges
coming out of $x$ and we may raise the effective conductance by identifying
all vertices other than $x$.

Fix some $a\in\R $ and set 
$$
f(x):=f_a(x):=a\,x-x^{2+{4\over d-2}}
\,.
$$
Since $f$ is concave, $f(M_n)$ is an $\F_n$-supermartingale.
Set $X_n:=1$ if $n=n_r$ for some $r\in\N$ and $\neg\ev A_n$ holds.
Otherwise, set $X_n:=0$.
We claim that $Y_n:=f(M_n)+b\,\sum_{j=0}^{n-1}X_j$ is an $\F_n$-supermartingale,
where $b>0$ is a certain constant depending only on $d$.
This will be established once
we show that
$$
\E\bigl[f(M_{n_r+1})\bigm|\F_{n_r}\bigr] \le f(M_{n_r})-b
\qquad\hbox{ on the event }
\neg \ev A_{n_r}\,.
\label e.gain
$$

Set $m:=n_r$ and assume that $\neg \ev A_{m}$ holds.
Note that since $\neg\ev A_m$ holds, we have $M_{m+1} \ge M_m$.
Set $L:=f(M_m)+f'(M_m)\,(M_{m+1}-M_m)$.
Then
$$
L-f(M_{m+1})=\int_{M_m}^{M_{m+1}} \bigl(f'(M_m)-f'(x)\bigr)\,dx
=\int_{M_m}^{M_{m+1}}\int_{M_m}^{x}-f''(y)\,dy\,dx\,.
\label e.Lf
$$
Since $f''(x)<0$, we get $L\ge f(M_{m+1})$.
Observe that there is a constant $C'=C'_d \ge 1$ such that $M_{n+1}\le C'\,M_n$.
(Given a unit flow from $S_{n+1}$ to $\infty$ in $\Z^d\setminus E_{n+1}$,
one can produce from it a unit flow from $S_n$ to $\infty$
in $\Z^d\setminus E_n$ by setting the flow appropriately
on $E_{n+1}\setminus E_n$. Since there is at most one edge in
$E_{n+1}\setminus E_n$ and the degrees in $\Z^d$ are bounded,
it is easy to see that the $\ell^2$ norm of the new flow is bounded
by some constant times the $\ell^2$ norm of the former flow.)
Therefore, $f''(M_m)/f''(x)$  is bounded on $[M_m,M_{m+1}]$ when $M_{m+1}\ge M_m$.
Recall that~\ref e.whener/ holds when $\te_r\in\fo$.
By our choice of $f$ and by~\ref e.Lf/, we therefore have
$$
 O_d(1)\,\bigl(L-f(M_{m+1})\bigr)\ge 1
\qquad\qquad\hbox{when }\te_r\in\fo\,.
$$
By~\ref e.per/, we therefore get $O_d(1)\,\E\bigl[L-f(M_{m+1})\bigm| \F_m\bigr]\ge 1$
(when $m=n_r$ and $\neg \ev A_{n_r}$ holds).
Since $M_n$ is a martingale, $\E[L\mid \F_m]=f(M_m)$, by the choice of $L$.
This proves~\ref e.gain/.

Now fix some $\bar M> C'\, M_0$.
Let $\tilde n:= \inf\{n\st M_n\ge \bar M/C'\}$, with the usual convention that
 $\tilde n:=\infty$ if $\forall n\ M_n<\bar M/C'$.
We now choose the constant $a$ in the definition of $f$ so that $f\ge 0$ on $[0,\bar M]$
and $f=0$ at the endpoints of this interval, namely, $a:= {\bar M}^{{d+2\over d-2}}$.
Set $Z_n:=\sum_{0 \le j < n} X_j$, where $n\in\N\cup\{\infty\}$.
Then
$$
\E[ Z_{n\wedge \tilde n}] = b^{-1} \,\E[Y_{n\wedge \tilde n}]-
b^{-1}\,\E[f(M_{n\wedge \tilde n})]\,.
$$
By our choice of $a$ and $\tilde n$, $f(M_{n\wedge \tilde n})\ge 0$.
Since $Y_n$ is a supermartingale, we get
$$
\E[Z_{n\wedge \tilde n}] \le b^{-1} \,Y_0= b^{-1}\,f(M_0) \le O_d(1)\,
{\bar M}^{{d+2\over d-2}}.
$$
The monotone convergence theorem implies that
the same bound applies to $\E[Z_{\tilde n}]$.
Therefore, for every $t>0$, we have
$$
\eqalign{
\P[Z_\infty>t]
&\le \P[Z_{\tilde n}>t]+\P[\tilde n\ne\infty]
\cr &
\le t^{-1}\,\E[Z_{\tilde n}]+\P\bigl[\sup\{ M_n\st n\in\N\} \ge \bar M/C'\bigr]
\cr &
\le O_d(1)\, t^{-1}\,
{\bar M}^{{d+2\over d-2}}
+C'M_0/\bar M\,.
}
$$
(The final inequality uses, say, Doob's optional stopping theorem.)
We now choose $\bar M:=t^{(d-2)/(2d)}$ for large $t$ and get
$\P[Z_\infty>t]\le O_d(1)\,t^{(2-d)/(2d)}$.
Observe that $O_d(1)\,Z_\infty$ bounds the Euclidean
diameter of the component of $o$ in $\fo$.
Note that (by taking a limit along an exhaustion) \ref l.bi-rooted/
implies that the diameter of the connected component of $o$ in $\fo$ under
$\wsf_\bfz$ stochastically dominates
the diameter of the union of the finite connected components of
$\fo\backslash \{\bfz\}$ under
the $\wsf$.
This completes the proof.
\Qed

\bsection{Conductance Criterion for One End}{s.gen}

In this section we consider a network $(G,c)$, and prove a rather general sufficient condition
for the $\wsf$ on $(G,c)$ to have one end for every tree a.s.

We recall that for a set of edges $F\subset\edges=\edges(G)$, we write
$|F|_c$ for $\sum_{e\in F}c(e)$, and for a set of vertices $K\subset \vertex=\vertex(G)$,
we write $\zalpha(K)$ for the sum of $c(e)$ over all edges $e\in\dedge$ having at least
one endpoint in $K$.

\procl p.mostgen
Let $(G, c)$ be a transient connected network and let
$V_0\subset V_1\subset\cdots$ be finite sets of vertices
satisfying $\bigcup_{j=0}^\infty V_j=\vertex$,
where $\vertex=\vertex(G)$.
Suppose that 
$$
\inf\big\{ \EC(v,\infty; G\backslash V_n)\st n\in\N,\, v\in\vertex\setminus
V_n\big\}>0
\label e.vertcond
$$
and
$$
\lim_{t\to\infty}
\inf\big\{ \EC(K,\infty; G\backslash V_n)\st n\in\N,\,
K\subset\vertex\setminus V_n,\, K \hbox{ finite},\,
\zalpha(K) > t\big\}=\infty\,.
\label e.condcond
$$
Then $\wsf$-a.s.\ every tree has one end.
\endprocl

\proof
Let $H_n$ be the subgraph of $G$ spanned by $V_n$.
Let $o\in\vertex$. With no loss of generality,
we assume that $o\in V_0$ and that
$\bde V_n\subset H_{n+1}$ for each $n\in\N$ (since we may
take a subsequence
of the exhaustion).
As in \ref s.Zdcase/, let $\fo$ be a sample from
$\wsf_o$ on $G$ and let $S_n$ be the set of vertices of the
connected component of $\fo\cap H_n$ containing $o$.
Let $E_n:=\edges(H_n)$, and let 
$E_n'$ be the set of edges in $E_n$ that have at least
one endpoint in $S_n$.
By \ref l.martingale/,
we know that the effective conductance $M_n$ from $S_n$ to
$\infty$ in $G\backslash E_n'$ is a non-negative  martingale and therefore
bounded a.s.
(Since in the formulation of the lemma it is assumed that
every edge of $E_1$ has an endpoint in $S_0$, in order to
deduce the above statement one has to first go through a 
procedure of examining edges one-by-one, as we have done
in \ref s.Zdcase/, but with the graphs $H_n$
used in place of the balls $B_r$.)

Fix some $n\in\N$.
Let $Z_n$ be a set of vertices such that
$\vertex\setminus Z_n$ is finite and contains $S_n$ and
$\EC(S_n,Z_n;G\backslash E_n') \le 2 \,M_n$.
By the definition of $M_n$, there is such a $Z_n$.
For every vertex $v\in \vertex$, let $h(v)$ denote the probability
that the network random walk on $G\backslash E_n'$
started at $v$ hits $Z_n$ before hitting $S_n$.
Let $F^n_0$ denote the set of edges in $\edge\setminus E_n'$ that join
vertices in $S_n$ to
vertices in $U_0^n:=\{v\in\vertex\setminus S_n\st h(v)\in[0,1/2]\}$, and let
$F^n_1$ denote the set of edges in $\edge\setminus E_n'$ that join vertices
in $S_n$ to vertices in $\{v\in\vertex\setminus S_n\st h(v)>1/2\}$.
We claim that a.s.
$$
\sup_{n\in\N} |F^n_0\cup F^n_1|_c <\infty\,.
\label e.bdcond
$$
We start by estimating $|F^n_0|_c$.  Set $H(v):=\max\{ 0,\, 2\,h(v)-1\}$.
The Dirichlet energy $D(H)$ of $H$ is at most $4$ times the Dirichlet energy
of $h$, which is bounded by $2\,M_n$. Since $H=1$ on $Z_n$, it follows that
$\EC(H^{-1}(0),Z_n;G\backslash E_n')\le D(H)\le 8\,M_n$.
Since $U_0^n \subseteq H^{-1}(0)$ and $E_n'\subset E_n$, it follows that
$\EC(U_0^n,\infty;G\backslash V_n)\le 8\,M_n$.
Since $M_n$ is a.s.\ bounded, it follows that $\sup_n
\EC(U_0^n,\infty;G\backslash V_n)<\infty$ a.s.,
and by \ref e.condcond/,
$\sup_n |F^n_0|_c \le \sup_n \zalpha(U_0^n) < \infty$ a.s.
Since 
$$
2\,M_n\ge D(h) = \sum_e c(e)\, \bigl(h(\etail e)-h(\ehead e)\bigr)^2
\ge \sum_{e\in F^n_1} c(e)\,(1/2)^2\,,
$$
we get $|F^n_1|_c\le 8\, M_n$, and so \ref e.bdcond/ follows.

Next, suppose that $e_1,e_2,\dots,e_k$ are all the edges in
$\edge\setminus E_n'$
joining $S_n$ to $\vertex\setminus V_n$, that is,
$F^n_0\cup F^n_1 = \{e_1,e_2,\dots,e_k\}$.
Let $v_j$ denote the vertex of $e_j$ outside of $S_n$ for $j=1,\dots,k$.
We now show that for $j=1,\dots,k$,
$$
\P\bigl[e_j\notin \fo \bigm| e_1,\dots,e_{j-1}\notin\fo, S_n\bigr]
\ge \exp\bigl(-c(e_j)/a\bigr)\,,
\label e.cutedge
$$
where $a>0$ is the left-hand side of~\ref e.vertcond/.
Let $\tilde G$ denote the graph $G\backslash (E_n'\cup\{e_1,\dots,e_{j-1}\})$.
Then $1$ minus the left-hand side of~\ref e.cutedge/ is equal to
$$
{c(e_j) \over \EC(v_j,\{\infty\}\cup S_n; \tilde G)}
\le 
{c(e_j)\over \EC(v_j,\infty;G\backslash V_n)+c(e_j)}
\le {c(e_j)\over a+c(e_j)}\,,
$$
which implies~\ref e.cutedge/
(using the inequality $(1+x)^{-1} \geq e^{-x}$,
valid for $x > -1$).
Now,~\ref e.cutedge/ and induction on $j$ imply 
$$
\P\bigl[ S_{n+1}=S_n \bigm| S_n\bigr] \ge \exp\bigl( -|F^n_0\cup F^n_1|_c/a\bigr)\,.
$$
Since by~\ref e.bdcond/ the right-hand side is a.s.\ bounded away from zero,
it follows that a.s.\ there is some $n\in\N$ such that $S_{n+1}=S_n$.
Consequently, the connected component of $o$ in $\fo$ is finite a.s.,
and this completes the proof by \ref l.bi-rooted/.
\Qed

\bsection{Isoperimetric Profile and Transience}{s.profile}

Our next goal is to prove that the $\wsf$ has 
one end a.s.\ in networks with a ``reasonable isoperimetric profile",
but first, we must discuss the relationship of the profile
to transience, which is the subject
of this section.
Define
$$
\prfl(\gh, A, t) :=
\inf \big\{ |\bde K|_c \st A \subseteq K,\,K \hbox{ is finite and connected},
t \le \zalpha(K) 
\big\}\,.
$$
We also abbreviate $\prfl(\gh,t):=\prfl(\gh,\emptyset,t)$.
Write 
$$
\bdei K := \big\{(x, y)\in\edge \st x \in K,\, y \hbox{ belongs to an infinite
component of } \gh \backslash K \big\}
\,.
$$
Then
$$
\prfl(\gh, A, t) =
\inf \big\{ |\bdei K|_c \st A \subseteq K,\,K \hbox{ is finite and connected},
t \le \zalpha(K) 
\big\}\,.
\label e.altinf
$$

The following result refines \ref b.Thomassen:profile/ and is adapted from
a similar result of \ref b.HeSchramm/.
Our proof is simpler than that of Thomassen, but
we do not obtain his conclusion of the existence of a transient
subtree.

\procl t.HS
Let $A$ be a finite set of vertices in a network $\gh$ with
$\zalpha\big(\vertex(\gh)\big) = \infty$.
Let $\prfl(t) := \prfl(\gh, A, t)$.
Define $s_0 := |\bdei(A)|_c$ and $s_{k+1} := s_k + \prfl(s_k)/2$ recursively for
$k \ge 0$.
Then 
$$
\ER(A, \infty) \le \sum_{k \ge 0} {2 \over \prfl(s_k)}
\,.
$$
\endprocl

This is an immediate consequence of the following analogue for finite 
networks.

\procl l.HSfinite
Let $a$ and $z$ be two distinct vertices in a finite connected network $\gh$.
Define
$$
\prfl(t) :=
\min \big\{ |\bde K|_c \st a \in K,\ z \notin K,\ \hbox{ $K$ is connected, }
t \le \zalpha(K) \big\}
$$
for $t \le \zalpha\big(\vertex(\gh)\setminus \{z\}\big)$ and 
$\prfl(t) 
:= \infty$ for $t > \zalpha\big(\vertex(\gh)\setminus \{z\}\big)$.
Define $s_0 := \sumC a$ and $s_{k+1} := s_k + \prfl(s_k)/2$ recursively for
$k \ge 0$.
Then 
$$
\ER(a, z) \le \sum_{k = 0}^{\infty} {2 \over \prfl(s_k)}
\,.
$$
\endprocl

\proof
Let $g:\vertex(G)\to\R$ be the function that is harmonic in $\vertex(G)\setminus\{a,z\}$
and satisfies $g(a)=0$, $g(z)=\ER(a,z)$.
Recall from~\ref s.back/ that $\nabla g$ is a unit flow from $a$ to $z$.
(To connect with electrical network terminology, note
that $-\nabla g$ is the unit current from $z$ to $a$ and $g$ is its voltage.)
For $t\ge 0$, let $W(t):=\bigl\{x\in\vertex\st g(x)\le t\bigr\}$,
and for $t'>t\ge 0$ let $E(t,t')$ be the set of directed edges 
from $W(t)$ to $\bigl\{x\in\vertex\st g(x)\ge t'\bigr\}$.
Define $t_0:=0$ and inductively,
$$
t_{k+1}:= \sup\bigl\{t\ge t_k \st |E(t_k,t)|_c\ge  |\bde W(t_k)|_c/2\bigr\}.
$$
Set $\bar k:=\min\{j\st z\in W(t_j)\}=\min\{j\st t_{j+1}=\infty\}$.
Fix some $k<\bar k$. Note that $\nabla g(e)\ge0$ for every $e\in\bde W(t_k)$
(where edges in $\bde W(t_k)$ are oriented away from $W(t_k)$).
Now
$$
\eqalign{
1 &= \sum_{e\in\bde W(t_k)}\nabla g(e) \ge
\sum_{e\in E(t_k,t_{k+1})} c(e)\,\bigl(g(\ehead e)-g(\etail e)\bigr)
\cr &
 \ge
\sum_{e\in E(t_k,t_{k+1})} c(e)\,(t_{k+1}-t_k)
\ge
(t_{k+1}-t_k)\,
{|\bde W(t_k)|_c\over 2}\,,
}
$$
where the last inequality follows from the definition of $t_{k+1}$.

Thus
$$
t_{k+1}-t_k\le 2/ \prfl\bigl(\zalpha(W_k)\bigr)\,,
\label e.prwk
$$
where we abbreviate $W_k:= W(t_k)$.
Clearly,
$$
\eqaln{
\zalpha(W_{k+1})
& =
\zalpha(W_k)+\zalpha(W_{k+1}\setminus W_k)
\ge
\zalpha(W_k)+|\bde W_k|_c-\sup_{t>t_{k+1}} |E(t_k,t)|_c
\cr &\ge
\zalpha(W_k)+{1\over 2}\,|\bde W_k|_c
\ge \zalpha(W_k) +{1\over 2}\, \prfl\big(\zalpha(W_k)\big)\,.
}
$$
Since $\prfl$ is a non-decreasing function, it follows by induction that
$\zalpha(W_k) \ge s_k$ for $k<\bar k$ and \ref e.prwk/ gives
$$
\ER(a, z)
=
g(z)
=
t_{\bar k}-t_0
\le
\sum_{k=0}^{\bar k-1} {2 \over \prfl\big(\zalpha(W_k)\big)}
\le
\sum_{k=0}^{\bar k-1} {2 \over \prfl(s_k)}
\,.
\Qed
$$

In the setting of \ref t.HS/,
it is commonly the case that $\prfl(t) = \prfl(\gh, A, t) \ge f(t)$ for some
increasing function $f$ on 
$\big[\pi(A), \infty\big)$
that satisfies $0 < f(t)
\le t$ and $f(2\,t) \le \alpha f(t)$ for some constant $\alpha$.
In this case, define $t_0 := \pi(A)$ and $t_{k+1} := t_k + f(t_k)/2$
recursively.
We have that $s_k \ge t_k$ and $t_k \le t_{k+1} \le 2\, t_k$, whence
for $t_k \le t \le t_{k+1}$, we have $f(t) \le f(2\, t_k) \le \alpha f(t_k)$,
so that
$$\eqalign{
\int_{\zalpha(A)}^\infty {4\,\alpha^2 \over f(t)^2} \,dt
&=
\sum_{k \ge 0} \int_{t_k}^{t_{k+1}} {4\,\alpha^2 \over f(t)^2} \,dt
\ge
\sum_{k \ge 0} \int_{t_k}^{t_{k+1}} {4 \over f(t_k)^2} \,dt
=
\sum_{k \ge 0} {4(t_{k+1} - t_k) \over f(t_k)^2}
\cr&=
\sum_{k \ge 0} {2 f(t_k) \over f(t_k)^2}
\ge
\sum_{k \ge 0} {2 \over \prfl(t_k)}
\ge
\sum_{k \ge 0} {2 \over \prfl(s_k)}
\ge
\ER(A, \infty)
\,.
}\label e.intbd
$$
This bound on the effective resistance is usually easier to estimate than
the one of \ref t.HS/.

We shall need the following fact, which states that a good isoperimetric
profile is inherited by some exhaustion.

\procl l.GoodSubset
Let $(\gh, c)$ be a connected locally finite
network such that 
$\lim_{t \to\infty} \prfl(\gh, \bp, t) = \infty$
for some fixed $\bp\in\vertex$ and
such that every infinite connected subset $K\subset\vertex(\gh)$ 
satisfies $\zalpha(K) = \infty$.
Then the network
$\gh$ has an exhaustion $\Seq{\gh_n}$ by finite connected subgraphs
such that 
$$
|\bde U \setminus \bde \vertex(\gh_n)|_c 
\ge
|\bde U|_c/2
\label e.minl
$$
for all $n$ and all finite $U \subset \vertex(\gh) \setminus
\vertex(\gh_n)$ and
$$
\prfl\big(\gh\backslash \vertex(\gh_n), t\big) 
\ge
\prfl(\gh, t)/2
\label e.GoodSubset
$$
for all $n$ and $t > 0$.
\endprocl

\proof
Given a finite connected $K \subset \vertex(\gh)$
containing $o$,
 let $W(K)$ minimize $|\bde L|_c$
over all finite sets $L\subset\vertex(\gh)$
that contain $K \cup \bdv K$
(here $\bdv K$ denotes the vertices outside of
$K$ neighboring some vertex in $K$); such a set $W(K)$
exists by our two hypotheses on $(G, c)$.
Moreover, $W(K)$ is connected since $K$ is.
Let $\gh' := \gh\backslash W(K)$ and write $\bde'$ for the edge-boundary
operator in $\gh'$.
If $U$ is a finite subset of vertices in $\gh'$, then $|\bde' U|_c \ge
|\bde U|_c/2$,
since if not, we would have $|\bde (W(K) \cup U)|_c<|\bde W(K)|_c$,
which contradicts the definition of $W(K)$.
Thus, $\prfl(\gh', t) \ge \prfl(\gh, t)/2$ for all $t > 0$.
It follows that we may define an exhaustion $G_n$ as the subgraphs
induced by a sequence $K_n$ defined 
recursively by $K_1 := W\big(\{\bp\}\big)$ and $K_{n+1} :=
W(K_n)$.
\Qed

\bsection{Isoperimetric Criterion for One End}{s.1end}

We now state and prove our general condition for the $\wsf$ trees to
have one end a.s.  After the proof, the range of its applicability will be discussed.

\procl t.profile1end
Suppose that $\gh$ is an infinite connected locally finite network.
Let $\prfl(t) := \prfl(\gh, t)$.
Suppose that $s_0:=\inf_{s>0}\prfl(s)>0$
and
$$
\sum_{k \ge 0} {1 \over \prfl(s_k)} < \infty\,,
\label e.isocond
$$
 where
$s_k$ is defined recursively by
$s_{k+1} := s_{k} + \prfl(s_{k})/2$ for $k\in\N$.
Then $\wsf$-a.s.\ every tree has only one end.
\endprocl

\proof
The proof will be based on an appeal to \ref p.mostgen/.
Note that the hypothesis \ref e.isocond/ certainly guarantees
$\lim_{t\to\infty}\prfl(\gh,\bp, t) \ge
\lim_{t\to\infty}\prfl(\gh, t)=\infty$ for every $\bp \in \vertex$.
Also, $s_0 > 0$ guarantees that $\zalpha(K) = \infty$ for $|K| = \infty$.
Thus, there exists an exhaustion $\Seq{\gh_n}$ satisfying 
the conclusion of \ref l.GoodSubset/, 
specifically, \ref e.GoodSubset/.
For $n \geq 0$, define $L_n := G \backslash \V(G_n)$. 
Consider some nonempty finite $A\subset\vertex(L_n)$.
 Define $r_0 := |\bd_{\edge(L_n)}^\infty(A)|_c$ and
$r_{k+1} := r_k + \prfl(L_n, A, r_k)/2$ recursively.
By~\ref e.GoodSubset/, $r_{k+1}\ge r_k+s_0/4$.
Therefore, we get in particular $r_4\ge s_0$.
We now prove by induction that if
$r_k\ge s_m$ for some $k,m\in\N$, then
$r_{k+2\ell}\ge s_{m+\ell}$ for every
$\ell\in\N$.
Indeed, by \ref e.GoodSubset/, we have $\prfl(L_n, r) \ge \prfl(\gh,
r)/2$ for all $r > 0$.
Apply this to $r := r_{k+2\ell+1}, r_{k+2\ell}$ in turn to obtain
$$\eqaln{
r_{k+2\ell+2}
&\ge
r_{k+2\ell+1} + \prfl(\gh, r_{k+2\ell+1})/4
\cr&\ge
r_{k+2\ell} + \prfl(\gh, r_{k+2\ell})/4
+ \prfl(\gh, r_{k+2\ell+1})/4
\cr&\ge
s_{m+\ell} + \prfl(\gh, s_{m+\ell})/2
=
s_{m+\ell+1}
}$$
if $r_{k+2\ell} \ge s_{m+\ell}$, which completes the induction step. 
The above results in particular give
$r_{k}\ge s_\ell$ when $k\ge 4+2\,\ell$.

Given $\epsilon > 0$, choose $m$ so that
$\sum_{k=m}^\infty 8 / \prfl(\gh, s_k) < \epsilon$
and choose $t$ so that $\prfl(\gh, t) \ge 2 s_m$.
Given any $n$ and any finite $A \subset \vertex(L_n)$ with $\zalpha(A) >
t$, we claim that $\ER(A, \infty; L_n) < \epsilon$, thereby establishing
\ref e.condcond/.
To see this, note that by \ref e.altinf/ and
 $\prfl(L_n, r) \ge \prfl(\gh, r)/2$ we have
$$
r_0
=
|\bd_{\edge(L_n)}^\infty(A)|_c
\ge
\prfl(L_n, t)
\ge
\prfl(\gh, t)/2
\ge
s_m
\,.
$$
Therefore, \ref t.HS/ and the above yield that
$$ 
\ER(A, \infty; L_n) \le
\sum_{k=0}^\infty {2 \over \prfl(L_n, r_k)} \le
\sum_{k=0}^\infty {4 \over \prfl(\gh, r_k)} \le
\sum_{k=m}^\infty {8 \over \prfl(\gh, s_k)} <
\epsilon\,. 
$$

A similar argument implies~\ref e.vertcond/ and completes the proof:
Let $A\subset\vertex(L_n)$ be a singleton.
Since we have $r_4\ge s_0$, we obtain
$$
\ER(A, \infty; L_n) \le
\sum_{k=0}^\infty {2 \over \prfl(L_n, r_k)} \le
\sum_{k=0}^\infty {4 \over \prfl(\gh, r_k)} \le
{16 \over s_0} + \sum_{k=0}^\infty {8 \over \prfl(\gh, s_k)} \,. 
\Qed
$$
\medskip

Which networks satisfy the hypothesis of~\ref t.profile1end/? Of course,
all non-amenable transitive networks do.
(By definition, a network $\gh$ is non-amenable 
if $\inf_{t>0}\prfl(\gh,t)/t>0$. It is transitive if its automorphism group
acts transitively on its set of vertices. It is quasi-transitive if the
vertex set breaks up into only finitely many orbits under
the action of the automorphism group.)
Although it is not obvious, so do all quasi-transitive transient graphs.
To show this,
we begin with the following slight extension of a result
due to \ref b.CoulSC/ and \ref
b.Sal:iso/.
Define the {\bf internal vertex boundary} of
a set $K$ as $\bdvi K := \{x \in K \st \texists {y \notin K} y \sim x\}$.

A locally compact group is called {\bf unimodular} if its
left Haar measure is also right invariant.
We call a graph $\gh$ {\bf unimodular} if its automorphism group
$\Aut(\gh)$ is unimodular, where
$\Aut(\gh)$ is given the weak topology generated by its action on $\gh$.
See \BLPSgip\ for more details on unimodular graphs.

\procl l.trnsIso
Let $\gh$ be an infinite
unimodular transitive graph. Let $\rho(m)$ be the smallest
radius of a ball in $\gh$ that contains at least $m$ vertices.
Then for all finite $K \subset \vertex$, we have 
$$
{|\bdvi K| \over |K|} \ge {1 \over 2\, \rho(2\,|K|)}
\,.
$$
\endprocl

\proof
Fix a finite set $K$ and let $\rho := \rho(2\,|K|)$.
Let $B'(x, r)$ be the ball of radius $r$ about $x$ excluding $x$ itself
and let $b := |B'(x, \rho)|$.
For $x, y, z \in \vertex(\gh)$, 
define $f_k(x,y,z)$ as the proportion of shortest paths 
from $x$ to $z$ whose $k$th vertex is $y$.
Let $S(x,r)$ be the sphere of radius $r$ about $x$.
Write $q_r := |S(x,r)|$.
Let $F_{r, k}(x, y) := \sum_{z \in S(x, r)} f_k(x, y, z)$.
Clearly, $\sum_y F_{r, k}(x,y)=q_r$ for every $x \in \vertex(\gh)$ and $r
\ge 1$.
Since $F_{r, k}$ is invariant under the diagonal action of the
automorphism group of $\gh$, the Mass-Transport Principle (\BLPSgip) gives
$\sum_x F_{r, k}(x,y)=q_r$ for every $y \in \vertex(\gh)$ and $r \ge 1$.
Now we consider the sum 
$$
Z_r :=
\sum_{x\in K} \;\sum_{z\in S(x, r) \setminus K}\; \sum_{y \in \bdvi K}\;
\sum_{k=0}^{r-1} f_k(x,y,z)
\,.
$$
If we fix $x\in K$ and $z\in S(x, r) \setminus K$,
then the inner double sum is at least 1, since if we fix any shortest path
from $x$ to $z$, it must pass through $\bdvi K$.
It follows that 
$$
Z_r \ge \sum_{x \in K} |S(x, r) \setminus K|
\,,
$$
whence, by the definitions of $\rho$ and $b$,
$$
Z := \sum_{r=1}^\rho Z_r
\ge
\sum_{x \in K} |B'(x, \rho) \setminus K|
\ge
\sum_{x \in K} |B'(x, \rho)|/2
=
|K| b/2
\,.
$$
On the other hand, if we do the summation in another order, we find 
$$\eqaln{
Z_r 
&=
\sum_{y \in \bdvi K}\;
\sum_{k=0}^{r-1}\; \sum_{x\in K}\; \sum_{z\in S(x, r) \setminus K} f_k(x,y,z)
\cr&\le
\sum_{y \in \bdvi K}\;
\sum_{k=0}^{r-1} \; \sum_{x\in \vertex(\gh)} \; \sum_{z\in S(x, r)}
f_k(x,y,z)
\cr&=
\sum_{y \in \bdvi K}\;
\sum_{k=0}^{r-1} \; \sum_{x\in \vertex(\gh)} F_{r, k}(x, y)
\cr&=
\sum_{y \in \bdvi K} \sum_{k=0}^{r-1} q_r
=
|\bdvi K|\, r\, q_r
\,.
}$$
Therefore, 
$$
Z \le \sum_{r=1}^\rho |\bdvi K|\, r \,q_r
\le |\bdvi K| \,\rho \, b
\,.
$$
Comparing these upper and lower bounds for $Z$, we get the desired result.
\Qed

An immediate consequence is the following bound:

\procl c.VolIso
If $\gh$ is a connected quasi-transitive graph
with balls of radius $n$ having at least $c\,
n^3$ vertices for some constant $c>0$, then 
$$
\prfl(\gh, t) \ge c' t^{2/3}
$$
for some constant $c'>0$ and all $t \ge 1$.
\endprocl

\proof
First assume that $\gh$ is transitive.
If $\gh$ is also amenable, then it is unimodular by \ref b.SoardiWoess:rws/.
Thus, the inequality follows from \ref l.trnsIso/.
If $\gh$ is not amenable, then the inequality is trivial by definition.

Now, suppose that $\gh$ is only quasi-transitive. Pick some vertex
$\bp\in\vertex(\gh)$, and let $\vertex'$ denote
the orbit of $\bp$ under the automorphism group of $\gh$. 
Let $r\in\N$ be such that every vertex in $\gh$ is within distance
$r$ of some vertex in $\vertex'$.
Let $\gh'$ be the graph on $\vertex'$ where two vertices are
adjacent if and only if the distance between them in $\gh$ is
at most $2\,r+1$. It is easy to verify that $\gh'$ satisfies
the assumptions of the corollary and is also transitive.
Consequently, we have $\prfl(\gh',t)\ge c''\,t^{2/3}$ for
some $c''>0$. The result now easily follows for $\gh$ as well.
\Qed

Since all quasi-transitive transient graphs have at least cubic volume
growth by a theorem of \ref b.Gromov:poly/ and \ref b.Trofimov/,
we may use \ref c.VolIso/, \ref e.intbd/ and \ref t.profile1end/  to obtain:

\procl t.wsfoneend If $G$ is a transient quasi-transitive network or is a
non-amenable network with $\inf_{x\in\vertex(\gh)}\sumC x >0$, then $\wsf$-a.s.\ every tree has only one end. 
\endprocl

In particular, we arrive at the following results that extend Theorem 12.7 of
\BLPSusf:

\procl t.H^dSF
Suppose that $G$ is a bounded-degree graph that is roughly isometric
to $\HH^d$ for some $d\ge 2$.
Then the $\wsf$ of $G$ has infinitely many trees a.s.,
each having one end a.s.
If $d=2$ and $G$ is planar,
then the $\fsf$ of $G$ has one tree with infinitely many ends a.s.
\endprocl 

Recall that when $d>2$ in the above setting, we have $\fsf=\wsf$.
(This follows from Theorems~7.3 and~12.6 in \BLPSusf.)

\proof
Rough isometry preserves non-amenability
(\ref b.Woess:book/, Theorem 4.7)
and $\HH^d$ is non-amenable when $d\ge 2$. Hence, $G$
is also non-amenable,
and \ref t.wsfoneend/ implies that the $\wsf$ has one end per tree a.s.
The fact that the $\wsf$ has infinitely many trees a.s.\ 
follows from Theorem~9.4 of \BLPSusf\ and the exponential
decay of the return probabilities
for random walks on non-amenable bounded-degree graphs.
(Stronger results are proved in Theorems 13.1 and 13.7 of \BLPSusf.)

When $d=2$ and $G$ is planar,
it is not hard to see that the planar dual of $G$ also has bounded degree.
This implies that it is roughly isometric to $G$, hence to $\HH^2$.
Thus, the above conclusions apply also to the $\wsf$ of the
dual of $G$.
Therefore, the claims about the $\fsf$ in the planar setting follow from Proposition~12.5 in \BLPSusf,
which relates the properties of the $\wsf$ on the dual to the properties
of the $\fsf$ on $G$.
\Qed

Finally, we conclude with some new questions that arise from our results.

\procl q.r-i1end
If $\gh$ and $\gh'$ are roughly isometric graphs and the wired spanning
forest in $\gh$ has only one end in each tree a.s., then is the same true
in $\gh'$?
\endprocl

\procl q.1end0-1
Is the probability that each tree has only one end equal to either 0 or 1
for both $\wsf$ and $\fsf$?
\endprocl

Consider the subgraph $\gh$ of $\Z^6$ spanned by the vertices
$$
\big(\Z^5\times\{0\}\big)\cup\big(\{(0,0,0,0,0),(2,0,0,0,0)\}\times \N\big)\,.
$$
This graph is obtained from $\Z^5$ by adjoining two copies of $\N$.
Let $\fo$ denote a sample from the $\wsf$ on $\gh$.
With positive probability, $x:=(0,0,0,0,0,0)$ and $y:=(2,0,0,0,0,0)$ are in the
same component of $\fo$. In that case, a.s.\ that component has $3$ ends while
all other trees in $\fo$ have one end. Also with positive probability,
$x$ and $y$ are in two distinct components of $\fo$. In that case,
a.s.\ each of these components have two ends while all other components
have one end. Thus, in particular, there are graphs such that the existence
of a tree in the $\wsf$ with precisely two ends has probability
in $(0,1)$.
This example naturally leads to the following question.

\procl q.excess
Define the {\bf excess} of a tree as the number of ends minus one.
Is the sum of the excesses of the components of the $\wsf$ equal to some
constant a.s.?
\endprocl

Note that the total number of ends of all trees is tail measurable, hence
an a.s.\ constant by Theorem 8.3 of \BLPSusf.
Since the number of trees is also a.s.\ constant by Theorem 9.4 of
\BLPSusf, it follows that \ref q.excess/ has a positive answer when the
number of trees is finite.

\procl q.anchored
Does our main result, \ref t.profile1end/, hold when $\prfl(\gh, t)$ is
replaced by $\prfl(\gh, \bp, t)$ for some fixed basepoint $\bp$? In
particular, does it hold when $\inf_t \prfl(\gh, \bp, t)/t > 0$? If so,
this would provide a new proof that Question 15.4 of \BLPSusf\ has a
positive answer in the case of bounded-degree Galton-Watson trees (by 
Corollary 1.3 of \ref b.ChPP:anchored/). This case of Galton-Watson trees
has already been established by \ref b.AL:urn/.
\endprocl

\procl q.recurrent-case
Suppose that $G$ is a quasi-transitive recurrent graph that is not roughly
isometric to $\Z$. Then \BLPSusf\ proves that the uniform spanning tree of
$G$ has only one end a.s. That proof is rather long; can it be simplified
and the result generalized?
\endprocl

\def\noop#1{\relax}
\input \jobname.bbl

\filbreak
\begingroup
\eightpoint\sc
\parindent=0pt\baselineskip=10pt

\def\emailwww#1#2{\par\qquad {\tt #1}\par\qquad {\tt #2}\smallskip}
\def\onlywww#1{\par\qquad {\tt #1}\smallskip}

Department of Mathematics,
Indiana University,
Bloomington, IN 47405-5701
\emailwww{rdlyons@indiana.edu}
{http://mypage.iu.edu/\string~rdlyons/}

Department of Mathematics,
University of California,
Davis, CA 95616
\emailwww{morris@math.ucdavis.edu}
{http://www.math.ucdavis.edu/\string~morris/}

Microsoft Research,
One Microsoft Way,
Redmond, WA 98052
\onlywww{http://research.microsoft.com/\string~schramm/}

\endgroup

\bye